\documentclass{article}
\usepackage{amsmath}
\usepackage{amsfonts}
\usepackage{amsthm}
\usepackage[all]{xy}
\newcommand{\CC}{\mathbb{C}}

\newcommand{\ZZ}{\mathbb{Z}}
\newtheorem{theorem}{Theorem}
\newtheorem*{theorem3}{Theorem 3}
\newtheorem{lemma}{Lemma}
\newtheorem{prop}{Proposition}
\newtheorem{cor}{Corollary}
\theoremstyle{definition}
\newtheorem{defn}{Definition}
\newtheorem*{defn2}{Definitions}
\theoremstyle{remark}
\newtheorem*{acn}{Acnowledgements}
\DeclareMathOperator{\ch}{ch}
\DeclareMathOperator{\pr}{pr}
\DeclareMathOperator{\inv}{inv}
\DeclareMathOperator{\Inv}{Inv}
\DeclareMathOperator{\Sym}{Sym}
\DeclareMathOperator{\dout}{out}
\DeclareMathOperator{\din}{in}
\DeclareMathOperator{\res}{res}
\DeclareMathOperator{\mult}{mult}
\begin{document}

\title{Tensor decompositions for $SL(2)$ and outerplanar graphs}

\author{Aleksandrs Mihailovs\\
Department of Mathematics\\
University of Pennsylvania\\
Philadelphia, PA 19104-6395\\
mihailov@math.upenn.edu\\
http://www.math.upenn.edu/$\sim$mihailov/
}
\date{\today}
\maketitle

\begin{abstract}
The main result of this article is the decomposition of tensor products 
of representations of $SL(2)$ in the sum of irreducible representations 
parametrized by outerplanar graphs. 
An {\em outerplanar graph} is a graph with the vertices $0,1,2, \dots, m$, 
edges of which can be drawn in the upper half-plane without intersections. 
We allow for a graph to have multiple edges, but don't allow loops. 
In more detail, 
\[
\rho_{d_1}\otimes\dots\otimes\rho_{d_m}=\underset{G}{\oplus} \thickspace T_G \thinspace,
\]
where $\rho_d$ denotes the irreducible representation of dimension $d+1$, and 
the direct sum is taken over all the outerplanar graphs of 
degrees $d_0,d_1,\dots,d_m$ with all possible values of $d_0$. $T_G$ is an 
irreducible subrepresentation of the type $\rho_{d_0}$, and we determine explicit 
formulas for the basis in the space of $T_G$ as well.\end{abstract}    
\setlength{\baselineskip}{1.5\baselineskip}

\section{Introduction}
In the classical quantum mechanics, each particle corresponds to an irreducible 
representations of $SL(2)$ of dimension $2s+1$, where $s \in \frac{1}{2} 
\thinspace \ZZ \thinspace$ is the spin 
of the particle. For example, a photon having spin 0 corresponds to the trivial 
representation of dimension $1$, and an electron having spin $1/2$ corresponds to the 
standard 2-dimensional representation of SL(2):
\begin{equation} \begin{pmatrix}a&b\\c&d\end{pmatrix} 
\begin{pmatrix}A\\ B\end{pmatrix} = \begin{pmatrix}aA+bB\\ cA+dB\end{pmatrix} . 
\label{1}
\end{equation}
The projections of the spin correspond to the elements of the fixed basis of the 
space of the representation, described below in section \ref{sec2} just before (\ref{9}).  
To study the systems of a few particles, one needs to know explicit formulas for the bases of 
the isotypic components of the tensor product of corresponding representations, and 
especially 
the basis of the isotypic component of the type of the trivial representation, i.\ e.\ 
the basis of the subspace of invariants of the tensor product of representations. 
These bases are of interest in the valency theory and a few other branches of 
physics as well as mathematics.   
   
The basis of the subspace of $SL(2)$-invariants of 
$V^{\otimes m}=V \otimes \dots \otimes V$ ($m$ times), 
where $V$ is the two-dimensional linear space 
with the standard action (\ref{1}) of $SL(2)$, was described in the terms of the 
$0$-outerplanar graphs in the classic work \cite{W} of  
Rumer, Teller and Weyl. This theory was developed and 
applied to the percolation theory by Temperley and Lieb \cite{TL},
to the knots theory and invariants of 3-manifolds by Jones \cite{J}, 
Kauffman \cite{Kauf}, Kauffman and Lins \cite{KL}, Wenzl \cite{Wen}, 
Jaeger \cite{Ja}, Lickorish \cite{L}, Masbaum and Vogel \cite{MV}, and others, 
to the quantum theory by Penrose \cite{P} and Moussouris \cite{Mou}, 
to quantum groups and the quantum link theory by Reshetikhin and Turaev \cite{RT},
Ohtsuki and Yamada \cite{OY}, Carter, Flath and Saito \cite{CFS} and others, 
especially to the theory of Lusztig's canonical bases 
\cite{Lu} by Khovanov and Frenkel \cite{FK}, Varchenko \cite{Var} and Frenkel, Varchenko and
Kirillov, Jr. 
\cite{FVK}. I am sorry that I am not able to cite everybody who made a contribution, 
because the literature on this topic is very extensive. 
I became familiar with the description of the basis of the invariants of the tensor products 
of any irreducible representations of $SL(2)$ in the terms of the $0$-outerplanar graphs from 
Kuperberg's work \cite{K}. Some people, and Frenkel and Khovanov \cite{FK} in particular, 
constructed bases not only in the invariants, but in the other components as well. The bases 
constructed in this work, are different.
In all the previous works, the proofs used straightening algorithm. 
They showed that each invariant could be expressed as a linear 
combination of invariants corresponding to the $0$-outerplanar graphs. 
That, after the calculation 
of the dimensions of an invariant space, and finding the number of $0$-outerplanar graphs, 
confirmed that the given invariants form a basis. 

In section \ref{sec2} 
I give new proof of the result of Rumer, Teller and Weyl \cite{W}.
The method of my proof is the following. 
I show that the invariants corresponding to the $0$-outerplanar graphs, 
are linearly independent. After the calculations of dimensions and numbers of 
$0$-outerplanar graphs, it proves that the given invariants form a basis. 

In section \ref{sec3} I give new formulation for the extension of Rumer, Teller and Weyl's 
work to the description of invariants of tensor products of any irreducible representations 
of $SL(2)$, using slightly different $0$-outerplanar graphs, than in \cite{K}. 
Then I give the new proof of the result, following 
the method described in the previous paragraph. 
Also, I give a few new explicit formulas of those invariants.

In section \ref{sec4} I give a few examples and a few new formulas of the invariants.

In section \ref{sec5} I introduce bases of all the isotypic components of tensor products of 
any irreducible representations of $SL(2)$, parametrized surprizingly by outerplanar graphs 
as well. Also, I describe the subdivisions of those bases, giving the decomposition of 
the isotypic component in the sum of the irreducible representations. The proofs are 
close to ones given in the previous sections.

Recall that all classes of equivalence of the irreducible polynomial finite dimensional 
representations of SL(2) are parametrized by nonnegative integers, and we can 
choose the natural actions $\rho_k$ of $SL(2)$ in the symmetric 
powers $S^kV, k=0, 1, 2, \dots$ as their representatives. 

\begin{defn2}
An {\em outerplanar graph} is a graph with the set of vertices $\{ 0,1,2, \dots, m 
\} \subseteq {\mathcal H}=\{ z \in \CC \: | \: \Im z \geq 0 \}$, edges of which can be drawn 
in the upper half-plane $\mathcal H$ without intersections. 
We allow for a graph to have multiple edges, but don't allow loops. 
Denote ${\mathcal O}(G)$ 
the set of directed graphs, underlying undirected graph of which is G. 
Let $G$ be an outerplanar graph and $g \in {\mathcal O}(G)$. For each vertix $k$ of $g$ 
denote 
\begin{equation}
x_k(g)=x^{d^{\dout}_k(g)}y^{d^{\din}_k(g)} \in S^{d_k}V , 
\label{i1}
\end{equation}
where $d^{\dout}_k(g)$ is the number of arrows in $g$, beginning in $k$, 
$d^{\din}_k(g)$ is the number of arrows, ending in $k$, and 
\begin{equation}
d_k=d^{\dout}_k(g)+d^{\din}_k(g) 
\label{i2}
\end{equation} 
is the degree of $k$. 
In other words, we put $x$ at the beginning of each arrow of $g$, 
and $y$ ---at the end, and multiply those $x$s and $y$s in each vertix. 
Denote 
\begin{equation}
b_g=x_1(g) \otimes \dots \otimes x_m(g) \in S^{d_1}V \otimes \dots 
\otimes S^{d_m}V 
\label{i3} 
\end{equation}
and for any nonnegative integer $i \leq d_0(G)$ denote 
\begin{equation}
t_{G,i}=\sum_{\substack{g\in {\mathcal O}(G)\\ d_0^{\din}(g)=i}}(-1)^{\inv g}b_g ,
\label{i4}
\end{equation}
where $\inv g$ is the number 
of inversions in $g$, i.\ e.\ , the number of arrows $(k,l)$ in $g$ with $k>l$. 
   
\end{defn2}

\begin{theorem3}
For any fixed $d_0, d_1, \dots, d_m$, tensors $t_{G,i}$ parametrized by all outerplanar graphs
with degrees $d_0, d_1, \dots, d_m$ and nonnegative integers $i \leq d_0$, 
form a basis in the isotypic component of the type 
$\rho_{d_0}$ in the representation $\rho_{d_1}\otimes \dots \otimes \rho_{d_m}$. For any 
fixed outerplanar graph $G$, the subspace $T_G$ spanned by the basis 
$(t_{G,0}, \dots, t_{G,d_0})$, are invariant; the linear homomorphism  
\begin{equation} \label{i5} 
s_G: S^{d_0}V \rightarrow T_G,\quad
x^iy^{d_0-i}\mapsto \frac{(-1)^i}{\binom{d_0}{i}} 
t_{G,i} 
\end{equation} 
defines the isomorphism of $\rho_{d_0}$ and the subrepresentation of $SL(2)$ in $T_G$, and
\begin{equation} \label{i6} 
S^{d_1}V\otimes\dots\otimes S^{d_m}V=\underset{G}{\oplus} \thickspace T_G \thinspace,
\end{equation} 
where the direct sum in the right hand side is taken over all the outerplanar graphs of 
degrees $d_0,d_1,\dots,d_m$ with all possible values of $d_0$.
\end{theorem3} 

\section{Invariants of tensor powers of\\ the standard represesentation}   
\label{sec2}  
Let $f$ be a field of characteristic $0$, and $SL(2)$ ---the group of 
$2 \times 2$ $f$-matrices with determinant 1, acting on $2$-dimensional 
linear $f$-space $V$ with basis $(x, y)$ by the standard way (\ref{1}): 
\begin{equation}
\begin{pmatrix} a&b\\ c&d \end{pmatrix}  (Ax+By)=
(aA+bB)x+(cA+dB)y , 
\label{2}
\end{equation}
i.\ e.\ 
\begin{equation} 
x \mapsto ax+cy ,\quad y \mapsto bx+dy . 
\label{3}
\end{equation}
Then $V \otimes V$ is $4$-dimensional linear $f$-space with the basis 
$(x \otimes x, x \otimes y, y \otimes x, y \otimes y)$, and $SL(2)$ acts 
on $V \otimes V$ through the tensor product of the standard actions (\ref{3}), 
i.\ e.\ 
\begin{equation} x \otimes x \mapsto (ax+cy) \otimes (ax+cy),
\label{4}
\end{equation}
and so on: 
\begin{equation}
\begin{split}
\begin{pmatrix} 
a&b\\ c&d \end{pmatrix}  
(Ax \otimes x + Bx \otimes y + Cy \otimes x +Dy \otimes y)&=
(a^2A+ab(B+C)+b^2D)x \otimes x\\ &+(acA+adB+bcC+bdD)x \otimes y\\
&+(acA+bcB+adC+bdD)y \otimes x\\  &+(c^2A+cd(B+C)+d^2D)y \otimes y. 
\end{split}
\label{5}
\end{equation} 

\begin{lemma}
The subspace of $SL(2)$-invariants 
of $V \otimes V$ is one-dimensional, and we can choose 
\begin{equation}
x \wedge y=x \otimes y - y \otimes x
\label{7}
\end{equation} 
as a basis element in that space.
\end{lemma}
\begin{proof}
From (\ref{5}), 
\begin{equation}
\begin{pmatrix} a&b\\ c&d \end{pmatrix} 
(x \otimes y - y \otimes x)=(ad-bc)(x \otimes y - y \otimes x)=
x \otimes y - y \otimes x.
\label{6}
\end{equation}
It means that $x \wedge y$ is invariant. 

Recall some fundamental facts about the 
representations of $SL(2)$. The word {\em representation} will mean below 
a polynomial finite dimensional linear representation over $f$.
Every representation 
of $SL(2)$ is equivalent to a sum of irreducible representations. 
All classes of equivalence of the irreducible 
representations are parametrized by nonnegative integers, and we can 
choose the natural actions $\rho_k$ of $SL(2)$ in the symmetric 
powers $S^kV, k=0, 1, 2, \dots$ as their representatives. 
In other words, $\rho_0$ is a trivial $1$-dimensional representation, 
$\rho_1$ is the standard representation (\ref{2}), and for $k>0$ the 
representation $\rho_k$ acts in the $(k+1)$-dimensional
linear $f$-space of the homogeneous polynomials of degree $k$ of two 
variables $x$ and $y$, which we can provide with the basis of monomials 
$(x^k, x^{k-1}y, \dots , y^k)$, such that 
\begin{equation}
x^k \mapsto (ax+cy)^k,\quad 
x^{k-1}y \mapsto (ax+cy)^{k-1}(bx+dy),\quad \dots ,\quad y^k \mapsto (bx+dy)^k . 
\label{9}
\end{equation}

By definition,
\begin{equation}
S^2V=(V \otimes V)/(fx \wedge y) .
\label{10}
\end{equation} 
It means that 
\begin{equation}
\rho_1 \otimes \rho_1 \simeq \rho_2 \oplus \rho_0 .
\label{11}
\end{equation}
As we see, the subspace of invariants is $1$-dimensional. 
\end{proof}

To study decompositions of tensor products of representations, it is convenient 
to consider characters, i.\ e.\ traces of the representations. For a diagonal 
matrix  $\bigl( \begin{smallmatrix} q&0\\ 0&q^{-1} \end{smallmatrix} 
\bigr)$ formulas (\ref{9}) give us 
\begin{equation}
x^k \mapsto q^kx^k, \quad 
x^{k-1}y \mapsto q^{k-2}x^{k-1}y, \quad \dots, \quad y^k \mapsto q^{-k}y^k .
\label{12}
\end{equation}
It means that
\begin{equation} 
\ch \rho_k=q^k+q^{k-2}+\dots+q^{-k}=\frac{q^{k+1}-q^{-(k+1)}}{q-q^{-1}}. 
\label{13}
\end{equation}

{\em Characters} are elements of the $\ZZ/2$-graded ring 
\begin{equation}
K=\ZZ [q+q^{-1}]=K^{\text{odd}} \oplus K^{\text{even}} . 
\label{14}
\end{equation}
Two representations are equivalent iff they have the same characters. 
For any representations $\sigma, \tau$ we have 
\begin{equation}
\ch \sigma \oplus \tau =\ch \sigma + \ch \tau 
\label{15}
\end{equation}
and 
\begin{equation}
\ch \sigma \otimes \tau = \ch \sigma \cdot \ch \tau .
\label{16}
\end{equation}
Characters of 
the representations $\rho_k$ with odd (even) $k$ form a basis of $\ZZ$-module 
$K^{\text{odd}}$ (or $K^{\text{even}}$, correspondingly), 
and it is useful to know the exact 
formulas for the coefficients in that basis of arbitrary odd (or even) character: 
\begin{equation}
\sum_{\text{odd (even)}\ k}C_kq^k=
\sum_{\text{odd (even)}\ k \geq 0}(C_k-C_{k+2}) \ch \rho_k . 
\label{17}
\end{equation}

\begin{lemma}The dimension of the subspace of $SL(2)$-invariants in $V^{\otimes m}$ 
equals to the Catalan number
\begin{equation}
c_n=\frac{(2n)!}{n!\thinspace (n+1)!}=\binom{2n}{n}-\binom{2n}{n-1}
\label{8}
\end{equation} 
for $m=2n$, and $0$ otherwise. 
\end{lemma}
\begin{proof}

Using formula (\ref{17}), we get
\begin{equation}
(\ch \rho_1)^m=(q+q^{-1})^m=\sum_{k \equiv m \bmod 2} \binom{m}{(m-k)/2} 
q^k= \sum_{k=0}^m c_m^{(k)} \ch \rho_k ,
\label{18}
\end{equation} 
where 
\begin{equation} 
c_m^{(k)}= \begin{cases}\binom{m}{(m-k)/2}-\binom{m}{(m-k)/2-1} & 
\text{if $m-k$ is even},\\ 0 & \text{otherwise}. \end{cases} 
\label{19}
\end{equation}
In particular, the dimension of the subspace of invariants in $V^{\otimes m}$
is equal to $c_m^{(0)}=c_n$ for $m=2n$, and 0 otherwise.
\end{proof} 

Catalan numbers $c_n$ count the number of 
$0$-outerplanar regular graphs of degree 1 with $m=2n$ vertices. 

\begin{defn}
A $0$-{\em outerplanar graph} is a graph with the set of vertices $\{ 1,2, \dots, m 
\} \subseteq {\mathcal H}=\{ z \in \CC \: | \: \Im z \geq 0 \}$, edges of which can be drawn 
in the upper half-plane $\mathcal H$ without intersections.  
\end{defn}  

For example, here are two of $c_3=5$ $0$-outerplanar regular graphs of degree 1 
with 6 vertices:

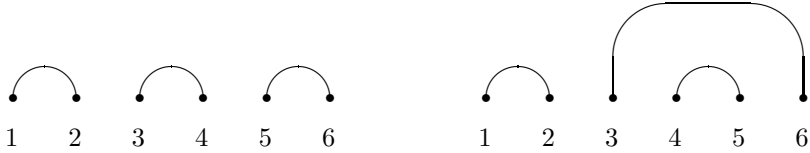
\begin{picture}(410, 54) 
  
\put(33,18){\circle*{3}}
\put(57,18){\circle*{3}}
\put(81,18){\circle*{3}}
\put(105,18){\circle*{3}}
\put(129,18){\circle*{3}}
\put(153,18){\circle*{3}}
\put(45,18){\oval(24,24)[t]}
\put(93,18){\oval(24,24)[t]}
\put(141,18){\oval(24,24)[t]}
\put(30,0){1}
\put(54,0){2}
\put(78,0){3}
\put(102,0){4}
\put(126,0){5}
\put(150,0){6}

\put(212,18){\circle*{3}}
\put(236,18){\circle*{3}}
\put(260,18){\circle*{3}}
\put(284,18){\circle*{3}}
\put(308,18){\circle*{3}}
\put(332,18){\circle*{3}}
\put(224,18){\oval(24,24)[t]}
\put(296,18){\oval(24,24)[t]}
\put(296,18){\oval(72,72)[t]}
\put(209,0){1}
\put(233,0){2}
\put(257,0){3}
\put(281,0){4}
\put(305,0){5}
\put(329,0){6}

\end{picture}

The tensor products of invariants are invariants. Edges of each $0$-outerplanar regular 
graph of degree $1$ give us $n$ pairs of numbers from $1$ to $m=2n$. If we take the
tensor product of the known invariants $x \wedge y$ in each of those pairs, we 
get an invariant in $V^{\otimes m}$. The count of such invariants is equal to 
the count of the considered $0$-outerplanar graphs, i.\ e.\ with the Catalan number $c_n$. 
And we know that the dimension of the subspace of invariants in $V^{\otimes m}$ 
equals the same Catalan number $c_n$. It is natural to suppose that the constructed 
invariants form a basis of the space of invariants. To prove that, it is enough 
to show that they are linearly independent. 

The standard basis $B$ of $V^{\otimes m}$ consists of $2^m$ tensor products 
$x_1 \otimes \dots \otimes x_m$ with $x_1, \dots, x_m \in \{ x, y \}$. 
We suppose that $B$ is ordered lexicographically. Denote ${\mathcal O}(G)$ 
the set of directed graphs, underlying undirected graph of which is G.
Let $G$ be an $0$-outerplanar regular graph of degree $1$, and $g \in {\mathcal O}(G)$.  
Denote 
\begin{equation}
b_g=x_1(g) \otimes \dots \otimes x_m(g) \in B
\label{20}
\end{equation}
setting $x_i(g)=x$, 
$x_j(g)=y$ for each edge $(i, j)$ in $g$. In other words, we put $x$ at the 
beginning of each arrow of $g$, and $y$ ---at its end. We can write the tensor 
product of $x \wedge y$ corresponding to $G$ that was introduced in the last paragraph  
as 
\begin{equation}
t_G=\sum_{g \in {\mathcal O}(G)}(-1)^{\inv g}b_g ,
\label{21}
\end{equation}
where $\inv g$ is the number 
of inversions in $g$, i.\ e.\ , the number of arrows $(i,j)$ in $g$ with $i>j$. 

\begin{theorem} Tensors $t_G$ parametrized by all regular $0$-outerplanar graphs
of degree $1$, form a basis in the subspace of invariants in the 
representation $\rho_1^{\otimes m}$.
\end{theorem} 
\begin{proof}
Notice that for each $G$ exists exactly one $g_0 \in {\mathcal O}(G)$ without inversions---with 
the orientation of each edge from the left to the right. Changing the orientation of the  
arrows of $g_0$ increases $b_g$ in the lexicografical order of $B$. It means that 
$b_{g_0}$ is the minimal element with a non-zero coefficient in the decomposition 
of $t_G$ in the basis $B$. Denote $b_G=b_{g_0}$. For the element of type $b_G$, 
we can reconstruct $G$, associating the left bracket with each $x$, the right bracket 
with each $y$, and connecting the corresponding left and right brackets. So, we have 
$c_n$ elements $b_G$ ---one for each $G$.

To prove the linear independence of the set of $t_G$, we can show that the rank 
of the $c_n \times 2^m$ matrix of the coefficients of $t_G$ in the basis $B$ is equal 
to $c_n$. To do that, we can find a non-zero $c_n \times c_n$ minor of that matrix. 
Consider the $c_n \times c_n$ submatrix with rows numbered by $G$ ordered the same way 
as $b_G$, and columns corresponding to $b_G$. As we noticed in the previous paragraph, 
$b_G$ is the first element with a nonzero coefficient in the row $G$, and  
this coefficient equals 1 by definition. So, this matrix is unipotent, its determinant 
equals $1$, that completes the proof of the linear independence of $t_G$.
\end{proof}

\section{General tensor invariants}
\label{sec3}
Now we allow for a graph to have multiple edges, but don't allow loops. 
Let $G$ be an $0$-outerplanar graph and $g \in {\mathcal O}(G)$. For each vertix $i$ of $g$ 
denote 
\begin{equation}
x_i(g)=x^{d^{\dout}_i(g)}y^{d^{\din}_i(g)} \in S^{d_i}V , 
\label{22}
\end{equation}
where $d^{\dout}_i(g)$ is the number of arrows in $g$, beginning in $i$, 
$d^{\din}_i(g)$ is the number of arrows, ending in $i$, and 
\begin{equation}
d_i=d^{\dout}_i(g)+d^{\din}_i(g) 
\label{23}
\end{equation} 
is the degree of $i$. 
In other words, we put $x$ at the beginning of each arrow of $g$, and $y$ ---at
the end, and multiply those $x$s and $y$s in each vertix. 
Denote 
\begin{equation}
b_g=x_1(g) \otimes \dots \otimes x_m(g) \in S^{d_1}V \otimes \dots 
\otimes S^{d_m}V 
\label{24} 
\end{equation}
and 
\begin{equation}
t_G=\sum_{g \in {\mathcal O}(G)}(-1)^{\inv g}b_g . 
\label{25}
\end{equation}

\begin{theorem}
\label{thm2}
For any fixed $d_1, \dots, d_m$, tensors $t_G$ parametrized by all $0$-outerplanar graphs
with degrees $d_1, \dots, d_m$, form a basis in the subspace of invariants in the 
representation $\rho_{d_1}\otimes \dots \otimes \rho_{d_m}$.
\end{theorem} 
\begin{proof} First, we can compare the dimension of the space of invariants and 
the number of $0$-outer-planar graphs with the given degrees. Using characters, we get
\begin{equation} 
\ch \rho_{d_1}\otimes \dots \otimes \rho_{d_m}=\prod_{i=1}^m(q^{d_i}+
q^{d_i-2}+\dots +q^{-d_i})=\sum_{k \equiv |d| \bmod 2}C_d^{(k)}q^k= 
\sum_{k=0}^{|d|} c_d^{(k)} \ch \rho_k , 
\label{26}
\end{equation} 
where $d=(d_1, \dots, d_m)$, $\thinspace |d|=d_1+\dots+d_m$ and
\begin{equation} 
c_d^{(k)}= \begin{cases}C_d^{(k)}-C_d^{(k+2)} & 
\text{if $|d|-k$ is even},\\ 0 & \text{otherwise}. \end{cases} 
\label{27}
\end{equation}
Applying the Clebsh-Gordon formula to the last two items of the product, 
we get 
\begin{equation}
\ch \rho_{d_{m-1}}\cdot \ch \rho_{d_m}=\ch \rho_{d_{m-1}+d_m}+\ch 
\rho_{d_{m-1}+d_m-2}+\dots+\ch \rho_{|d_m-d_{m-1}|} . 
\label{28}
\end{equation}
It gives us the recursion relation
\begin{equation}
c_d^{(k)}=c_{(d_1,\dots,d_{m-2},d_{m-1}+d_m)}^{(k)}+c_{(d_1,\dots,d_{m-2},
d_{m-1}+d_{m-2}-2)}^{(k)}+\dots+c_{(d_1,\dots,d_{m-2},|d_m-d_{m-1}|)}^{(k)} 
\label{29} 
\end{equation} 
for any $k$, and in particular, for $k=0$ defining the dimension of the space of 
invariants. That relation makes it possible to decrease $m$, and for $m=1$ we have the
initial conditions 
\begin{equation}
c_{(d)}^{(k)}=\delta_{dk} , 
\label{30} 
\end{equation} 
because the representations $\rho_d$ and $\rho_k$ are irreducible.

Denote $c_d$ the number of $0$-outerplanar graphs with degrees of vertices $d=(d_1,...,d_m)$. 
For any $0$-outerplanar graph with the given degrees, 
contracting $m$ to $m-1$, and deleting after that all the loops at the last vertix, 
we get an $0$-outerplanar graph with the degrees $(d_1,\dots,d_{m-2},d_{m-1}+d_m-2a_{m-1,m})$, 
where $a_{m-1,m}$ denotes the number of edges between $m-1$ and $m$ ---each edge has two 
ends---it explains the coefficient $2$. The number $a_{m-1,m}$ can be any integer between $0$ 
and $\min \{ d_{m-1},d_m\}$. That means the last degree in the new $0$-outerplanar graph can 
be $d_{m-1}+d_m, d_{m-1}+d_m-2, \dots, |d_m-d_{m-1}|$. This procedure is reversible:
for any $0$-outerplanar graph with admissible new degrees we can construct an $0$-outerplanar 
graph with degrees $d$, moving the last $d_m$ ends of edges from $m-1$ to $m$ and adding 
$a_{m-1,m}$ edges between $m-1$ and $m$. I wrote 'the last'
bearing in mind the natural order on the ends of edges of a vertix $i$ of 
an $0$-outerplanar graph: we can suppose that a small half-circle with the center in $i$, is
lying in the upper half-plane $\mathcal H$, which intersects with each edge ending in $i$, 
exactly in one point, and we can order these points on the half-circle clockwise, and these 
ends of edges---correspondingly.  So, we have a bijection between the set of $0$-outerplanar 
graphs with degrees d, 
and the union of the sets of $0$-outerplanar graphs of the mentioned above degrees. 
It gives us the same recursion relation as (\ref{29}) for $c_d^{(k)}$:
\begin{equation}  
c_d=c_{(d_1,\dots,d_{m-2},d_{m-1}+d_m)}+c_{(d_1,\dots,d_{m-2},
d_{m-1}+d_{m-2}-2)}+\dots+c_{(d_1,\dots,d_{m-2},|d_m-d_{m-1}|)} . 
\label{31} 
\end{equation}
For $m=1$ we have only one $0$-outerplanar graph---with one vertix $1$ and without edges. 
So, we have the same initial conditions as (\ref{30}) for $c_{(d)}^{(0)}$: 
\begin{equation}
c_{(d)}=\delta_{d0} . 
\label{32}
\end{equation}
Thus, we proved that the dimension $c_d^{(0)}$ of the space of invariants is equal to the 
number $c_d$ of the $0$-outerplanar graphs with the given degrees $d$.

The proof of the linear independence of the set of $t_G$ is exactly the same as for the case 
of regular $0$-outerplanar graphs of degree $1$. Denote $B$ the standard basis of 
$S^{d_1}V \otimes \dots \otimes S^{d_m}V$, 
consisting of $(d_1+1)\cdots (d_m+1)$ tensor products 
$x_1 \otimes \dots \otimes x_m$ with $x_i \in \{ x^{d_i}, x^{d_i-1}y, \dots , y^{d_i} \}$ 
for $i=1, \dots , m$. We suppose that $B$ is ordered lexicographically. We have $b_g \in B$ 
for all $g \in {\mathcal O}(G)$. 
For each $G$ exists exactly one $g_0 \in {\mathcal O}(G)$ without inversions - 
with orientation each edge from the left to the right. Changing the orientation of the  
arrows of $g_0$ increases $b_g$ in the lexicografical order of $B$. It means that 
$b_{g_0}$ is the minimal element with a non-zero coefficient in the decomposition 
of $t_G$ in the basis $B$. Denote $b_G=b_{g_0}$. In any vertix $i$ of $g_0$ we have, 
first of all the incoming arrows, and then - all the arrows coming out, in the order of the 
ends of the edges described above. 
Thus, for an element of type 
\begin{equation}
b_G=x^{d^{\dout}_1}y^{d^{\din}_1} \otimes \dots \otimes x^{d^{\dout}_m}y^{d^{\din}_m} , 
\label{33}
\end{equation}  
we can reconstruct $G$, associating the sequence of $d^{\din}_i$ in the right brackets followed 
by $d^{\dout}_i$ in the left brackets, with each vertix $i$ , and the connecting corresponding 
left and right brackets. So, we have $c_d$ elements $b_G$ ---one for each $G$. 
Consider the $c_d \times c_d$ submatrix with rows numbered by $G$ ordered the same way 
as $b_G$, and columns corresponding to $b_G$, of the $c_d \times |B|$ matrix of the 
coefficients of $t_G$ in the basis $B$.  As we noticed, 
$b_G$ is the first element with a nonzero coefficient in the row $G$, and  
this coefficient equals 1 by definition. So, this matrix is unipotent, its determinant 
equals $1$, that completes the proof of the linear independence of $t_G$. 

To complete the proof of the theorem, we have to show that the tensors $t_G$ are 
$SL(2)$-invariant. The representation $\rho_{d_1} \otimes \dots \otimes \rho_{d_m}$ is a 
subrepresentation of $\rho^{\otimes m}$ in $(SV)^{\otimes m}$, where 
$SV=f[x,y]$ and $\rho$ denotes the standard changing of variables. $(SV)^{\otimes m}$ 
is an algebra, and $\rho^{\otimes m}$ commutates with its multiplication 
\begin{equation}
(a_1 \otimes \dots \otimes a_m)(b_1 \otimes \dots \otimes b_m)=a_1b_1 \otimes \dots \otimes 
a_mb_m , 
\label{34}
\end{equation}
because $\rho$ commutates with the multiplication in $SV$. In particular, a 
product of invariants is invariant. It follows from the definition that 
\begin{equation}
t_G=\prod_{(i,j) \in E(g_0)}\iota_{ij}(x \wedge y) , 
\label{35}
\end{equation}
where $E(g_0)$ denotes the set of the arrows of $g_0$, and 
\begin{equation}
\iota_{ij}: V \otimes V \rightarrow (SV)^{\otimes m}, \quad 
v \otimes w \mapsto 1^{\otimes (i-1)} \otimes v \otimes 1^{\otimes (j-i-1)} 
\otimes w \otimes 1^{\otimes (m-j)} , 
\label{36}
\end{equation}
i.\ e.\ it puts $v$ on the $i$-th place and $w$ on the 
$j$-th place of the tensor product. 
Now, $\iota_{ij}$ commutate with the actions of $SL(2)$, $x \wedge y$ is invariant, 
and $t_G$, the product of invariants, is invariant.   
\end{proof}

\begin{prop}
\label{prop1}
Let $E(G)$ be the set of the edges of the $0$-outerplanar graph $G$ with $m$ vertices. 
For any disjoint union 
\begin{equation}
E(G)=E_1\amalg \dots \amalg E_k 
\label{36a}
\end{equation}
we have 
\begin{equation} 
t_G=t_{G_1}\dots t_{G_k} ,
\label{36b}
\end{equation}
where $G_1, \dots, G_k$ are the $0$-outerplanar graphs with $m$ vertices, with the 
sets of the edges $E_1, \dots, E_k$ correspondingly.
\end{prop}
\begin{proof}
It follows from (\ref{35}) immediately.
\end{proof}

\section{Examples}
\label{sec4}
If $m=1$, we have exactly one $0$-outerplanar graph: with the vertix $1$ and without edges. 
It means that the subspace of invariants of $\rho_k$ is nontrivial iff $k=0$. 
In that case, the corresponding invariant is $1 \in f$.
\begin{cor}
\label{cor1} 
The subspace of invariants of $S^{d_1}V \otimes S^{d_2}V $ 
is nontrivial iff $d_1=d_2=a$ for an integer $a \geq 0$. In that 
case, the space of invariants is one-dimensional with the basis element 
\begin{equation} 
(x \wedge y)^a=\sum_{i=0}^a(-1)^i \binom{a}{i} x^{a-i}y^i \otimes x^iy^{a-i} \in 
S^aV \otimes S^aV . 
\label{37}
\end{equation}
\end{cor}
\begin{proof}
If $m=2$, the degrees $d_1, d_2$ of any $0$-outerplanar graph $G$ are equal.
For any integer $a$ there is a unique $0$-outerplanar graph with 
two vertices of degrees $d_1=d_2=a$.
Formula (\ref{37}) follows from (\ref{35}).
\end{proof}

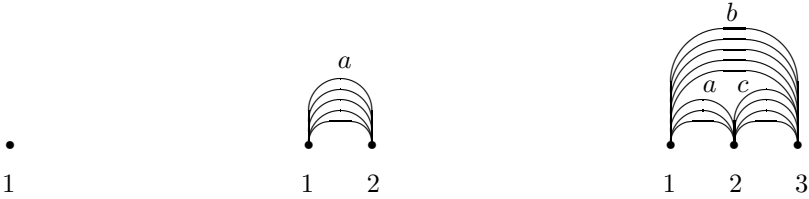
\begin{picture}(410, 80)
\put(30,18){\circle*{3}}
\put(27,0){1}

\put(143,18){\circle*{3}} 
\put(167,18){\circle*{3}}
\put(155,18){\oval(24,18)[t]}
\put(154,47){$a$}
\put(155,18){\oval(24,26)[t]}
\put(155,18){\oval(24,34)[t]}
\put(155,18){\oval(24,42)[t]}
\put(155,18){\oval(24,50)[t]}
\put(140,0){1}
\put(165,0){2}

\put(280,18){\circle*{3}}
\put(304,18){\circle*{3}}
\put(328,18){\circle*{3}}
\put(292,18){\oval(24,18)[t]}
\put(292,38){$a$}
\put(292,18){\oval(24,26)[t]}
\put(292,18){\oval(24,34)[t]}
\put(316,18){\oval(24,18)[t]}
\put(305,38){$c$}
\put(316,18){\oval(24,26)[t]}
\put(316,18){\oval(24,34)[t]}
\put(316,18){\oval(24,42)[t]}
\put(304,18){\oval(48,56)[t]}
\put(301,65){$b$}
\put(304,18){\oval(48,64)[t]}
\put(304,18){\oval(48,72)[t]}
\put(304,18){\oval(48,80)[t]}
\put(304,18){\oval(48,88)[t]}
\put(277,0){1}
\put(302,0){2}
\put(327,0){3}
\end{picture}

\begin{cor}
\label{cor2} 
The subspace of invariants of $S^{d_1}V \otimes S^{d_2}V \otimes 
S^{d_3}V$ is nontrivial iff the integers
$d_1, d_2$ and $d_3$ could be the lengths of the sides of a triangle 
(perhaps, degenerate, i.\ e.\ sides can have zero lenghth, or a vertix can be situated in 
the opposite side) with an even perimeter. In that 
case, the space of invariants is one-dimensional with the basis element 
\begin{equation} t_G=
((x \wedge y)^a\otimes 1)(x\otimes 1 \otimes y - y\otimes 1 \otimes x)^b(1\otimes (x 
\wedge y)^c) , 
\label{40}
\end{equation}
where
\begin{equation} 
a=\frac{d_1+d_2-d_3}{2}, \quad b=\frac{d_1+d_3-d_2}{2}, \quad c=\frac{d_2+d_3-d_1}{2} . 
\label{39}
\end{equation}
\end{cor}
\begin{proof} 
If $m=3$, we have a $0$-outerplanar graph with $a$ edges between 1 and 2, with $b$ edges 
between 1 and 3, and $c$ edges between 2 and 3 for any nonnegative integers $a$, $b$, and $c$. 
In this case 
\begin{equation}
d_1=a+b, \quad d_2=a+c, \quad d_3=b+c , 
\label{38} 
\end{equation}
and conversely (\ref{39}). Formula (\ref{40}) follows from (\ref{35}).
\end{proof}

In the degenerate cases some of $a, b, c$ equal 0. If all of them equal 0, we have a 
trivial representation. If two of them equal 0, say $b=c=0$, the situation is the same as 
for $m=2$:
\begin{equation}
t_G=(x \wedge y)^a \otimes 1=\sum_{i=0}^a(-1)^i \binom{a}{i} x^{a-i}y^i \otimes x^iy^{a-i} 
\otimes 1 \in S^aV \otimes S^aV \otimes f . 
\label{41}
\end{equation}
If one of them equals 0, say $b=0$, it follows from Corollary \ref{cor1} and 
Proposition \ref{prop1}, that 
\begin{equation} \begin{split} 
t_G&=((x \wedge y)^a\otimes 1)(1\otimes (x \wedge y)^c)\\ &=\sum_{i=0}^a\sum_{j=0}^c 
(-1)^{i+j} \binom{a}{i} \binom{c}{j} x^{a-i}y^i \otimes x^{c+i-j}y^{a-i+j} \otimes 
x^jy^{c-j} \in S^aV \otimes S^{a+c}V \otimes S^cV .
\end{split} \label{42}
\end{equation} 
It is a special case of the following 
\begin{prop}
\label{prop2}
Let $\Gamma_G$ denote the graph without multiple edges,
with the same vertices as $G$,  
in which an edge between two vertices exists 
iff there are edges between these vertices in $G$. Then  
\begin{equation}
t_G=\sum_{i_1=0}^{a_1}\dots \sum_{i_N=0}^{a_N} (-1)^{i_1+\dots+i_N} 
\binom{a_1}{i_1} \dots \binom{a_N}{i_N}b_{i_1,\dots, i_N} , 
\label{43}
\end{equation}
where $N$ is the number of edges of $\Gamma_G$, and for each edge $i_k$ of $\Gamma_G$ we 
denote $a_k$ the number of the edges of $G$ connecting the same vertices, and 
$b_{i_1,\dots, i_N}=b_g$ for any directed graph 
$g \in {\mathcal O}(G)$ with the exactly $i_1$ inverted 
arrows connecting the vertices of the 1-st edge of $\Gamma_G$, \dots, the exactly $i_N$ 
inverted arrows connecting the vertices of the $N$-th edge of $\Gamma_G$. 
If $\Gamma_G$ is a tree, i.\ e.\ a graph (not necessary connected) without cycles, 
then $(\ref{43})$ is the decomposition of $t_G$ in the elements of the basis $B$. 
\end{prop} 
\begin{proof}
It follows from Corollary \ref{cor1} and 
Proposition \ref{prop1}, that formula (\ref{43}) is true for any 
$0$-outerplanar graph $G$. We have to prove that for the considered case in which 
$\Gamma (G)$ is a tree, tensors $b_{i_1,\dots,i_N}$ for 
different $i_1, \dots, i_N$ are not equal. We'll use the induction on $N$. 
For $N=0$ we have only one item in the right hand side of (\ref{43}), so the 
statement is true. Let $N>0$ and we know that the statement is true for $N-1$. 
Choose the vertix $k$
of $\Gamma (G)$ with degree 1 (it exists because $\Gamma (G)$ is a tree). 
Then $k$-th component of $b_{i_1,\dots,i_N}$ equals by formula (\ref{22}) to  
\begin{equation}
x_k=
\begin{cases}x^{a_K-i_K}y^{i_K} &\text{if $k<l$,}\\
x^{i_K}y^{a_K-i_K} &\text{if $k>l$,}
\end{cases}
\label{44}
\end{equation} 
where $K=\{ k, l\} $ is the unique edge of $\Gamma (G)$ having an end $k$, and another its 
end is $l$. It means that $b_{i_1,\dots,i_N}$ are different for different $i_K$. And if 
$i_K$ are the same, then by induction, considering $G$ without edges connecting $k$ and 
$l$, we see that $b_{i_1,\dots,i_N}$ without $k$-th component, and without some factor 
on $l$-th place, are different for different $i_1, \dots, i_N$. They will differ  
on the same place after adding equal $k$-th components and multiplying $l$-th component on 
equal factors. 
\end{proof}

If $\Gamma (G)$ is not a tree,  some $b_{\dots}$ in formula (\ref{43}) can be the same. 
For instance, let $a=b=c=1$. Then 
\begin{equation} (-1)^{\inv g}b_g=
\begin{cases}
+x^2 \otimes xy \otimes y^2 &\text{if} \thickspace g= 
\xymatrix{1 \ar[r] \ar @/^/ [rr] & 2 \ar[r] & 3} ,\\
-xy \otimes xy \otimes xy &\text{if} \thickspace g= 
\xymatrix{1 \ar[r] & 2 \ar[r] & 3 \ar @/_/ [ll]} ,\\
-x^2 \otimes y^2 \otimes xy &\text{if} \thickspace g= 
\xymatrix{1 \ar[r] \ar @/^/ [rr] & 2 & 3 \ar[l]} ,\\ 
+xy \otimes y^2 \otimes x^2 &\text{if} \thickspace g= 
\xymatrix{1 \ar[r] & 2 & 3 \ar[l] \ar @/_/ [ll]} ,\\ 
-xy \otimes x^2 \otimes y^2 &\text{if} \thickspace g= 
\xymatrix{1 \ar @/^/ [rr] & 2 \ar[l] \ar[r] & 3} ,\\ 
+y^2 \otimes x^2 \otimes xy &\text{if} \thickspace g= 
\xymatrix{1 & 2 \ar[l] \ar[r] & 3 \ar @/_/ [ll]} ,\\ 
+xy \otimes xy \otimes xy &\text{if} \thickspace g= 
\xymatrix{1 \ar @/^/ [rr] & 2 \ar[l] & 3 \ar[l]} ,\\ 
-y^2 \otimes xy \otimes x^2 &\text{if} \thickspace g= 
\xymatrix{1 & 2 \ar[l] & 3 \ar[l] \ar @/_/ [ll]} . 
\end{cases}
\label{45}
\end{equation} 
The items corresponding to cycles, cancel, and $t_G$ is the sum of 
the other items. It is a special case of the following 
\begin{cor}
\label{cor3}
In the situation of the Corollary \ref{cor2} for $d_1=d_2=d_3=2a$,
\begin{equation}
t_G=(x^2 \wedge xy \wedge y^2)^a .
\label{46}
\end{equation}
\end{cor} 
\begin{proof} For $a=1$ see (\ref{45}), and for the other $a$ it follows from 
this special case and Proposition \ref{prop1}.
\end{proof}
 
\section{Tensor decompositions}
\label{sec5}

In this section we use slightly different outerplanar graphs than before: with 
vertices starting from $0$ instead of 1:

\begin{defn}
An {\em outerplanar graph} is a graph with the set of vertices $\{ 0,1,2, \dots, m 
\} \subseteq {\mathcal H}=\{ z \in \CC \: | \: \Im z \geq 0 \}$, edges of which can be drawn 
in the upper half-plane $\mathcal H$ without intersections.  
\end{defn} 

We still allow for a graph to have multiple edges, but don't allow loops. 
Also we'll use the same formulas (\ref{22}--\ref{24}) to define $x_i(g)$ and $b_g$ 
for an directed graph $g \in {\mathcal O}(G)$. Note that formula (\ref{24}) for $b_g$ 
doesn't include $x_0(g)$. For any nonnegative integer $i \leq d_0(G)$ denote 
\begin{equation}
t_{G,i}=\sum_{\substack{g\in {\mathcal O}(G)\\ d_0^{\din}(g)=i}}(-1)^{\inv g}b_g .
\label{47}
\end{equation}   

\begin{theorem}
\label{thm3}
For any fixed $d_0, d_1, \dots, d_m$, tensors $t_{G,i}$ parametrized by all outerplanar graphs
with degrees $d_0, d_1, \dots, d_m$ and nonnegative integers $i \leq d_0$, 
form a basis in the isotypic component of the type 
$\rho_{d_0}$ in the representation $\rho_{d_1}\otimes \dots \otimes \rho_{d_m}$. For any 
fixed outerplanar graph $G$, the subspace $T_G$ spanned by the basis 
$(t_{G,0}, \dots, t_{G,d_0})$, are invariant; the linear homomorphism  
\begin{equation} \label{48} 
s_G: S^{d_0}V \rightarrow T_G,\quad
x^iy^{d_0-i}\mapsto \frac{(-1)^i}{\binom{d_0}{i}} 
t_{G,i} 
\end{equation} 
defines the isomorphism of $\rho_{d_0}$ and the subrepresentation of $SL(2)$ in $T_G$, and
\begin{equation} \label{48a} 
S^{d_1}V\otimes\dots\otimes S^{d_m}V=\underset{G}{\oplus} \thickspace T_G \thinspace,
\end{equation} 
where the direct sum in the right hand side is taken over all the outerplanar graphs of 
degrees $d_0,d_1,\dots,d_m$ with all possible values of $d_0$.
\end{theorem} 
\begin{proof}
First we'll compare the  multiplicity of $\rho_{d_0}$ with the count of the outerplanar graphs.
By formulas (\ref{26}--\ref{27}), this multiplicity equals 
\begin{equation}
c_{d_1,\dots,d_m}^{(d_0)}= \begin{cases}C_{d_1,\dots,d_m}^{(d_0)}-C_{d_1,\dots,d_m}
^{(d_0+2)} & 
\text{if $d_0+d_1+\dots+d_m$ is even},\\ 0 & \text{otherwise}, \end{cases} 
\label{49}
\end{equation}
where the coefficients $C_{\cdot}^{\cdot}$ are defined by the generating function
\begin{equation} 
\sum_k C_{d_1,\dots,d_m}^{(k)}q^k= \prod_{i=1}^m(q^{d_i}+
q^{d_i-2}+\dots +q^{-d_i}) .
\label{50}
\end{equation} 
We can easily transform an outerplanar graph to a $0$-outerplanar graph and backwards, shifting
it to the right, or to the left, correspondingly. Thus the number of the outerplanar 
graphs with degrees $d_0,\dots,d_m$ is equal to the number of the $0$-outerplanar graphs with 
the same degrees, but shifted: the degree of $1$ must be $d_0$, \dots , the degree of $m+1$ 
must be $d_m$, i.\ e.\ equal to 
\begin{equation}
c_{d_0,d_1,\dots,d_m}=c_{d_0,d_1,\dots,d_m}^{(0)} .
\label{51}
\end{equation}
This equality was proved on the first step of the proof of Theorem \ref{thm2}. Now 
\begin{equation}
\begin{split}
c_{d_1,\dots,d_m}^{(d_0)}&=C_{d_1,\dots,d_m}^{(d_0)}-C_{d_1,\dots,d_m}^{(d_0+2)}= 
C_{d_1,\dots,d_m}^{(d_0)}-C_{d_1,\dots,d_m}^{(-d_0-2)}\\ 
&=\res_0 (q^{-d_0-1}-q^{d_0+1})\prod_{i=1}^m(q^{d_i}+q^{d_i-2}+\dots +q^{-d_i})\\ 
&=\res_0 (q^{-1}-q)\prod_{i=0}^m(q^{d_i}+q^{d_i-2}+\dots +q^{-d_i})\\
&=C_{d_0,d_1,\dots,d_m}^{(0)}-C_{d_0,d_1,\dots,d_m}^{(2)}=c_{d_0,d_1,\dots,d_m}^{0}=
c_{d_0,d_1,\dots,d_m} ,
\end{split}
\label{52}
\end{equation}
where $\res_0$ denotes a residue in 0. It means that the number of outerplanar graphs of 
degrees $d_0, d_1, \dots, d_m$ coincides with the multiplicity of $\rho_{d_0}$ in 
$\rho_{d_1}\otimes \dots \otimes \rho_{d_m}$. Notice, that using formula (\ref{51}) 
we are able to obtain that result very easily from the following equality as well:
\begin{equation}
\mult_{\rho}\tau=\dim \Inv(\rho^{\ast} \otimes \tau) 
\label{53}
\end{equation}
for any irreducible representation $\rho$ and representation $\tau$, where $\mult_{\rho}
\tau$ is the multiplicity of $\rho$ in $\tau$, $\rho^{\ast}$ denotes the representation 
conjugated to $\rho$, and $\dim \Inv$ is the dimension of the subspace of invariants. In our 
case 
\begin{equation}
c_{d_1,\dots,d_m}^{(d_0)}=\mult_{\rho_{d_0}}\rho_{d_1}\otimes\dots\otimes\rho_{d_m}=
\dim \Inv (\rho_{d_0}\otimes\rho_{d_1}\otimes\dots\otimes\rho_{d_m})=c_{d_0,d_1,\dots,d_m}
^{(0)} .
\label{53a}
\end{equation}

On the second step of the proof we'll show that the set of $t_{G,0}$ is linearly independent. 
For each $G$ exists exactly one $g_0 \in {\mathcal O}(G)$ without inversions - 
with orientation each edge from the left to the right. Changing the orientation of the  
arrows of $g_0$ with non-zero ends increases $b_g$ in the lexicografical order of $B$. 
It means that 
$b_{g_0}$ is the minimal element with a non-zero coefficient in the decomposition 
of $t_{G,0}$ in the basis $B$. Denote $b_G=b_{g_0}$. 
Note that for the $0$-outerplanar graph $G'$, obtained from $G$ by shifting it to the right on 
$1$, we have 
\begin{equation} \label{54} 
b_{G'}=x^{d_0}\otimes b_G .
\end{equation}
In the proof of Theorem \ref{thm2}, we checked that all $c_{d_0,d_1,\dots,d_m}$ elements 
$b_{G'}$ are different for all the $0$-outerplanar graphs of degrees $d_0, d_1, \dots, d_m$. 
Thus, by formula (\ref{54}), all $c_{d_0,d_1,\dots,d_m}$ elements $b_G$ are different for all 
the outerplanar graphs of degrees $d_0, d_1, \dots, d_m$.  
Consider the $c_{d_0,d_1,\dots,d_m}\times c_{d_0,d_1,\dots,d_m}$ submatrix with rows numbered
by $G$ ordered the same way 
as $b_G$, and columns corresponding to $b_G$, of the $c_{d_0,d_1,\dots,d_m}\times |B|$ 
matrix of the coefficients of $t_{G,0}$ in the basis $B$.  As we noticed, 
$b_G$ is the first element with a nonzero coefficient in the row $G$, and  
this coefficient equals 1 by definition. So, this matrix is unipotent, its determinant 
equals $1$, that completes the proof of the linear independence of $t_{G,0}$. 

Third step: $T_G$. First consider the case of a star with the center in $0$, i.\ e.\ 
$G$ having degrees $d_0=m,\thickspace d_1=\dots=d_m=1$:

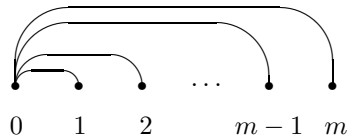
\begin{picture}(410,50) 
\put(147,18){\circle*{3}}
\put(171,18){\circle*{3}}
\put(195,18){\circle*{3}}
\put(210,18){\dots}
\put(243,18){\circle*{3}}
\put(267,18){\circle*{3}}
\put(159,18){\oval(24,12)[t]}
\put(171,18){\oval(48,24)[t]}
\put(195,18){\oval(96,48)[t]}
\put(207,18){\oval(120,60)[t]}
\put(145,0){0}
\put(169,0){1}
\put(194,0){2}
\put(230,0){$m-1$}
\put(264,0){$m$}
\end{picture}

In this case 
\begin{equation} \label{55}
t_{G,i}=(-1)^i \binom{m}{i} \Sym (x^{\otimes i}\otimes y^{\otimes (m-i)}) , 
\end{equation}
where $x^{\otimes i}=x\otimes\dots\otimes x$ ($m$ times) and $y^{\otimes (m-i)}=
y\otimes\dots\otimes y$ ($m-i$ times) analogously, and  
\begin{equation} \label{56}
\Sym (x_1\otimes\dots\otimes x_m)=\frac{1}{m!}\sum_{\sigma \in {\mathcal S}_m}
x_{\sigma(1)}\otimes\dots\otimes x_{\sigma(m)} ,
\end{equation} 
where ${\mathcal S}_m$ is the symmetric group of permutations of $1, \dots, m$.
Thus 
\begin{equation} \label{57}
T_G=S_mV=\Sym(V^{\otimes m}).
\end{equation}
By the definition of the actions of $SL(2)$ in $S^mV$ and $V^{\otimes m}$, the homomorphism 
\begin{equation} \label{58} 
s_G: S^mV \rightarrow S_mV, \quad x^iy^{m-i} \mapsto \Sym(x^{\otimes i}\otimes 
y^{\otimes (m-i)}) 
\end{equation}
is intertwining, therefore $s_G$ is an isomorphism and $T_G=S_mV$ is an invariant subspace. 

Now let $G$ be a star with the center in $0$ and multiple edges, i.\ e.\  
\begin{equation} \label{59}
d_0=d_1+\dots+d_m .
\end{equation}
Suppose for simplicity that all $d_1, \dots, d_m$ are not zero. 
Denote $\tilde{G}$ the covering $G$ star with the the same degree $d_0$ of $0$ but without 
multiple edges, i.\ e.\ with $d_1=\dots=d_{d_0}=1$. Consider a projection 
\begin{equation} \label{60}
\Pr=\pr_{d_1}\otimes\dots\otimes\pr_{d_m}: V^{d_0}\rightarrow S^{d_1}V\otimes\dots\otimes
S^{d_m} ,
\end{equation}
where 
\begin{equation} \label{61}
\pr_d: V^{\otimes d} \rightarrow S^dV 
\end{equation}
is the projection defining $S^d$. One can check that for any nonnegative integer $i\leq d_0$, 
\begin{equation} \label{62} 
t_{G,i}=\Pr(t_{\tilde{G},i}) .  
\end{equation}
Hense
\begin{equation} \label{63} 
T_G=\Pr(T_{\tilde{G}}) .
\end{equation}
By definitions (\ref{60}--\ref{61}), the projection $\Pr$ is intertwining with the 
$SL(2)$-actions. Therefore, $T_G$ is an invariant subspace since $T_{\tilde{G}}$ is 
an invariant subspace. The representation of $SL(2)$ in $T_G$ is a factor-representation 
of the irreducible representation of $SL(2)$ in $T_{\tilde{G}}$, and $T_G$ is non-zero 
space, because it contains $t_{G,0}$ that is not zero since the set of the tensors 
$t_{G,0}$ is linearly independent---it was stated on the first step of the proof of the 
theorem. It means that the representation of $SL(2)$ in $T_G$ is equivalent $\rho_{d_0}$, 
and 
\begin{equation} \label{64}
s_G=\Pr \circ \thinspace s_{\tilde{G}}
\end{equation} 
is the corresponding isomorphism. 
 
Adding a few vertices of degrees $0$ to $G$ doesn't change the situation, because of 
the canonical isomorphism: 
\begin{equation} \label{64a}
\alpha: S^{d_1}V\otimes\dots\otimes S^{d_m}V \rightarrow 1^{\otimes n_0}\otimes 
S^{d_1}V\otimes 1^{\otimes n_1}\otimes\dots\otimes 1^{\otimes n_{m-1}}\otimes S^{d_m}
\otimes 1^{\otimes{n_m}} , 
\end{equation}
adding 1's to every place corresponding to a vertix of degree $0$. In particular, 
we have 
\begin{equation} \label{64b}
t_{G_+\!,i}=\alpha (t_{G,i}), \quad T_{G_+}=\alpha (T_G) ,
\end{equation}
and 
\begin{equation} \label{64c}
s_{G_+}=\alpha \circ s_G ,
\end{equation}
where $G_+$ is an outerplanar graph obtained from $G$ by adding $n_0$ vertices of degree $0$ 
between $0$ and $1$, \dots, $n_k$ vertices of degree $0$ between $k$ and $k+1$, \dots, 
$d_m$ vertices of degree $0$ after $m$. 
 
Now, one can check that 
for an arbitrary outerplanar graph $G$, and for any nonnegative integer $i\leq d_0$, 
\begin{equation} \label{65}
t_{G,i}=t_{G_{\ast}\!,i} \thinspace t_{G_0} ,
\end{equation}
where $G_0$ is the $0$-outerplanar graph obtained from $G$ by deleting the vertex $0$ 
together with all the edges ending in $0$, and $G_{\ast}$ is the star obtained from $G$ 
by deleting all the edges of $G_0$. Since we stated in Theorem \ref{thm2} that $t_{G_0}$ is 
invariant, and the multiplication is intertwining, it follows from (\ref{65}) 
that the linear homomorphism 
\begin{equation} \label{66}
\mu : T_{G_{\ast}} \rightarrow T_G, \quad t \mapsto t \thinspace t_{G_0}
\end{equation}
is intertwining and surjective. Further, $T_G$ is nonzero, because it contains $t_{G,0}\neq 0$, 
see above. Since the representation of $SL(2)$ in $T_{G_{\ast}}$ is irreducible, and equivalent 
$\rho_{d_0}$, the representation of $SL(2)$ in $T_G$, equivalent its factor-representation, 
is equivalent $\rho_{d_0}$ as well, and 
\begin{equation} \label{67} 
s_G=\mu\circ s_{G_{\ast}}
\end{equation} 
is the corresponding isomorphism. 

Final step. We already know that all the subspaces $T_G$ are invariant, and the count of 
them is exactly the same as necessary for (\ref{48a}). Now we can use the following 
\begin{lemma} \label{lem3}
For any isotypic representation $T$ of $SL(2)$ and for any $q \in f$ that is not a root of 
unity, all the eigenspaces of $T\bigl( \begin{smallmatrix} q&0\\ 0&q^{-1} \end{smallmatrix} 
\bigr)$ have the same dimension equal to the number of the irreducible components of $T$.
\end{lemma}
\begin{proof}
For $\rho_d$, the basis elements $x^d,\dots,x^iy^{d-i},\dots,y^d$ are the eigenvectors of 
$\rho_d\bigl( \begin{smallmatrix} q&0\\ 0&q^{-1} \end{smallmatrix} \bigr)$ with the 
eigenvalues $q^d,\dots,q^{d-2i},\dots,q^{-d}$ that are all different if $q$ is not a root 
of unity. Therefore, for $T$ being a direct sum of $n$ representations equivalent to $\rho_d$, 
decomposing each component in the sum of the indicated above eigenlines (i.\ e.\ 
eigenspaces of dimension $1$), we get the 
decomposition of the space of $T$ in the sum of $n$ eigenlines with the eigenvalue $q^d$, 
$\dots$, $n$ eigenlines with the eigenvalue $q^{d-2i}$, $\dots$, $n$ eigenlines with the 
eigenvalue $q^{-d}$. It means that the dimension of each eigenspace is $n$. 
\end{proof} 

Consider the sum of $T_G$ for all $G$ with fixed degrees $d_0, d_1, \dots, d_m$. It is a space 
of isotypic representation of the type $\rho_{d_0}$. It follows from formula (\ref{48}) proved 
above, that each $t_{G,i}$ is an eigenvector with the eigenvalue $q^{d-2i}$ for the matrix 
in Lemma \ref{lem3}. In particular, the eigenspace with the eigenvalue $q^{-d}$ contains 
$c_{d_0,d_1,\dots,d_m}$-dimensional space spanned by $t_{G,0}$. By Lemma \ref{lem3}, the 
number of irreducible components in the considered space is not less than 
$c_{d_0,d_1,\dots,d_m}$. We know that the number of the irreducible components of the type 
$\rho_{d_0}$ in $\rho_{d_1}\otimes\dots\otimes\rho_{d_m}$ is $c_{d_0,d_1,\dots,d_m}$. 
It means that the considered sum of $c_{d_0,d_1,\dots,d_m}$ subspaces $T_G$ is direct, 
and it is the space of the isotypic component of the type $\rho_{d_0}$ in 
$\rho_{d_1}\otimes\dots\otimes\rho_{d_m}$. Taking the direct sum of the isotypic components, 
we get (\ref{48a}).
\end{proof}

\begin{prop}
\label{prop3}
Let $E(G)$ be the set of the edges of the outerplanar graph $G$ with $m$ vertices. 
For any disjoint union 
\begin{equation}
E(G)=E_1\amalg \dots \amalg E_k 
\label{67a}
\end{equation}
and for any nonnegative integer $i \leq d_0(G)$ we have 
\begin{equation} 
t_{G,i}=\sum_{i_1+\dots+i_k=i}t_{G_1,i_1}\dots t_{G_k,i_k} ,
\label{68}
\end{equation}
where $G_1, \dots, G_k$ are the $0$-outerplanar graphs with $m$ vertices, with the 
sets of the edges $E_1, \dots, E_k$ correspondingly, and we suppose in the sum (\ref{68}) 
that $i_1\leq d_0(G_1), \dots, i_k\leq d_0(G_k)$. 
\end{prop}
\begin{proof}
For a star without multiple edges, see (\ref{55}), it follows directly from the 
the definition (\ref{47}) of $t_{G,i}$. Using formula (\ref{65}) and 
Proposition \ref{prop1} we obtain the general case.
\end{proof}

\begin{acn}
I would like to thank Fan Chung Graham for her encouragement and  
series of very useful suggestions;  
Alexandre Kirillov for familiarizing me with Greg Kuperberg's work and giving me 
the opportunity to present my work at his seminar October 2, 1996; 
Greg Kuperberg himself for writing a nice article and giving me a lot of references; 
Murray Gerstenhaber and Michael Khovanov for useful discussions, 
and last 
but not least, my Gorgeous and Brilliant Wife, Bette, for her total support and love.
\end{acn}

\end{document}